\documentclass[11pt,noinfoline,addressatend,keywordsasfootnote]{article}

\RequirePackage[OT1]{fontenc}

\RequirePackage[amsthm,amsmath,natbib]{imsart}
\usepackage{times,amsfonts,amssymb}

\startlocaldefs

\newtheorem{theorem}{Theorem}[section]
\newtheorem{lemma}[theorem]{Lemma}
{\bf}{\rm}
\numberwithin{remarks}{section}
\newtheorem{proposition}[theorem]{Proposition}
\newtheorem{definition}{Definition}
\numberwithin{theorem}{section} \numberwithin{definition}{section}

\def\X{\mathbb{X}}

\def\Z{\mathbb{Z}}

\def\R{\mathbb{R}}
\def\G{\mathcal{G}}
\def\1{1}

\def\P#1{ \mathbb{P} \left( {#1} \right) }

\def\E#1{\mathbb{E}\left( {#1} \right) }

\endlocaldefs

\date{June 8, 2006}

\begin{document}

\begin{frontmatter}

\title{On the Critical Behavior at the Lower Phase Transition of the Contact Process}
\runtitle{On the Lower Phase Transition of the Contact Process}
\runauthor{M. Aizenman and P. Jung}

\author{\fnms{Michael} \snm{Aizenman}\ead[label=e1]{aizenman@princeton.edu}}
\thanks{Supported in part by NSF grant DMS 0602360}
\affiliation{Physics and Mathematics Departments, Princeton
University}
\address{Michael Aizenman \\
Departments of Physics and Mathematics\\
Princeton University, Princeton, NJ 08544 \\  \printead{e1}}
{\rm \and}

\author{\fnms{Paul}
 \snm{Jung}\ead[label=e2]{pjung@math.cornell.edu}}
 \affiliation{Mathematics Department, Cornell University}
\address{Paul Jung \\  591 Malott Hall\\
Cornell University\\ Ithaca, NY 14850\\ \printead{e2}}


\medskip
\medskip

\begin{abstract}
We present general results for the contact process by a method 
which applies to all transitive graphs of bounded degree, including 
graphs of exponential growth. The model's infection rates are varied 
through a control parameter, for which two natural transition points 
are defined as: i. $\lambda_T$, the value up to which the infection 
dies out exponentially fast if introduced at a single site, and 
ii. $\lambda_H$, the threshold for the existence of an invariant 
measure with a non-vanishing density of infected sites. It is 
shown here that for all transitive graphs the two thresholds coincide.  
The method, which proceeds through partial differential inequalities 
for the infection density, yields also generally valid bounds on two
related critical exponents.  The main results discussed here were
established by Bezuidenhout and Grimmett~\cite{BG2} in an extension to the continuous-time process of the discrete-time analysis of Aizenman and Barsky~\cite{AB}, and of the partially similar results of Menshikov~\cite{M}.
 The main novelty  here is in the direct derivation of the partial differential inequalities by an argument which is formulated for the continuum.
\end{abstract}

\begin{keyword}[class=AMS]
\kwd[Primary ]{60K35}
\end{keyword}
\begin{keyword}
\kwd{contact process} \kwd{phase transition} \kwd{critical
behavior} \kwd{interacting particle system} \kwd{oriented percolation}
\end{keyword}

\end{frontmatter}

\def\eq#1{eq.~(\ref{#1})}
\def\too#1{\mathop {\longrightarrow}_{#1}}
\maketitle
 \newpage
\nocontentsline
\tableofcontents

\newpage

\section{Introduction and statement of the main results}
Since its introduction by Harris~\cite{H}, the contact process has attracted interest as a model for the spread of ``infection". The model undergoes a phase transition which is reached by varying the ratio of the infection rate to the healing rate, which in our notation is $\lambda : 1$. The small $\lambda$ regime can be
characterized by the finiteness of the ``susceptibility",
$\chi(\lambda)$, which is the total time lost to infection within
the population  if an infection is introduced at a single site. A
duality argument allows to conclude that if $\chi(\lambda)<\infty$
then the infection dies out even if initially the entire
population was infected.   As $\lambda$ approaches the  edge of the
regime $\{\lambda :   \chi(\lambda)<\infty\}$  the susceptibility
diverges and the contact process exhibits critical behavior with
characteristics similar to those observed in models of statistical
mechanics.

Upon analysis, it turns out to be
a generally valid statement that right past the point of divergence of $\chi(\lambda)$
a homogeneous contact process enters the phase at which there is a stationary measure with persistent infection.  The technique presented below allows to establish this basic feature for the contact processes on the broad  class of {\em transitive graphs}. This is the class of graphs which are invariant under the action of a symmetry group which acts   transitively.    Included in the collection are some graphs for which  the contact process is known to exhibit more than one transition, in the sense explained below.
In this generality, the basic properties of the model
include the following:
\begin{enumerate}
\item For $\lambda$ with $\chi(\lambda) < \infty$  the probability
that infection from a single site will persist in the population
for time $t$ decays exponentially in time.   In models for which
the rates for the direct transmission have suitable exponential
decay, the probability that the infection would reach distance $d$
away also exhibits, in that phase, exponential decay in the
distance. \item At the edge of the above region
$\chi(\lambda)\nearrow\infty$, i.e., the model exhibits
criticality at the point
\begin{equation}
\lambda_T \  :=  \ \sup\{\lambda \, | \,  \chi(\lambda) < \infty
\} \, .
\end{equation}
which is named for Temperley. \item For $\lambda > \lambda _T $
the model is in the phase at which the infection persists. The
threshold for the latter condition has been recognized by the term
$\lambda_H $, for Hammersley.  Thus the above statement amounts to
the coincidence of the two points: $\lambda_T \ = \  \lambda_H $.
\item At the transition point, which can now be denoted simply
$\lambda_c$, the model exhibits critical behavior with
characteristic exponents which in general are bounded by their
`mean-field' values.  The bounds are realized in certain
situations.
\end{enumerate}

In regard to the spread of infection,
the contact process can be viewed as  oriented percolation.  That offers a helpful perspective, as the above characteristics are shared by transitive percolation models with or without orientation.
It was in that context that the characteristics 1.-3. of the phase diagram were initially established for the discrete-time version of the models, in two different and independently derived methods, presented in the works of  Menshikov~\cite{M} and Aizenman and Barsky~\cite{AB}.   The argument of \cite{M} was limited to graphs of subexponential growth, such as $\Z^d$.  The method of \cite{AB} readily extends to transitive graphs and
yields also additional information on the critical exponents, which follow through partial differential inequalities on which more is said below.  However, both analyses were initially presented for only the discrete-time version of the models.   The extension to the continuous-time contact process was accomplished in the work of Bezuidenhout and Grimmett~\cite{BG2}, through a detailed control of the (1D) continuum limit.
Our main goal here is to present a direct extension of the method of \cite{AB} to the contact process in terms which are natural for the continuum, casting the argument  in the generality of the transitive graphs of arbitrary  growth rate.

\subsection{The model and its parameters}

We shall introduce the model in the context of transitive graphs.
Before giving a formal description of the generator of the time evolution,
 let us set some notation.  For  $\mathcal{G}=(V,E)$ a
connected transitive graph with vertex set $V$ and edge collection
$E$, the contact process $ \{ A_t \}_{t\in\R}$ is a random
time-dependent collection of subsets of $V$ describing a set of
infected sites.  If the initial set of infected sites is given by
$B$ at time $T$, then the corresponding measure on $ \{ A_t
\}_{t\in\R}$ is given by $\mathbb{P}^{(B,T)}(\cdot)$. Let
\begin{equation*}
c(y,A)=\left\{
\begin{array}{ll}
1&\text{if }y\in A\\
\lambda\sum_{x\in A}J_{x,y}&\text{if }y\not\in A
\end{array}
\right..
\end{equation*}
where $J_{x,y}$ is a translation invariant kernel with
\begin{equation*}
|J|= \sum_y  J_{x,y} < \infty.
\end{equation*}
The generator of the contact process $A_t$ is formally given by
\begin{equation}\label{generator}
L f(A)=\sum_{x}{c(x,A)[f(A\bullet \{x\} )-f(A)]} \, ,
\end{equation}
where  $\bullet $ denotes the symmetric difference operation,
i.e., for $x\in V$: $A\bullet \{x\}=A\cup \{x\} $ when $x\not\in
A$ and $A\bullet \{x\}=A\backslash \{x\} $ when $x\in A$.

Often times the kernel is $1$ when $x$ and $y$ are neighbors and
$0$ otherwise; the finiteness of $|J|$ implies that the vertex
degree is finite. A simple description of the process in this case
is as follows.  At time $t$, the set of infected vertices is
denoted by $A_t$.  A vertex heals independently with exponential
rate $1$, while uninfected vertices become infected at exponential
rate $\lambda$ times the number of infected neighbors.

Two significant quantities which reflect properties of the model are: \\
1.  the infection density of the upper invariant measure
\begin{eqnarray} \theta_+(\lambda)  \ :=\
 \lim_{T\to - \infty}  \mathbb{P}^{(V,T)}( o \in A_0 )  \   =\
\lim_{T\to - \infty}  \mathbb{P}^{(\{o\},T)}( A_0\neq\emptyset  )
\end{eqnarray}
where the second equation is by duality and  stationarity,   \\
2.   the susceptibility
\begin{equation}
\chi(\lambda) \  := \int_0^\infty \mathbb{E}^{(\{o\},0) }(|A_t|)
dt<\infty \, ,
\end{equation}
which, by Fubini, equals the expected value of the sum of the
times lost to infection at the different sites.

The contact process exhibits a  number of different phases,
depending on the control parameter $\lambda$.  Some of the
thresholds of interest are defined as follows.
\begin{definition}
\begin{eqnarray*}
\lambda_T&:=&\sup\{\lambda:\, \chi(\lambda)<\infty \,  \} \\
\lambda_H&:=&\sup\{\lambda:\, \theta_+(\lambda)=0\, \}\\
\lambda_{GN}&:=&\sup\{\lambda : \, \mathbb{P}^{(\{o\},0)}( o \in
A_t) \too{t \to \infty} 0\, \}
\end{eqnarray*}
\end{definition}
Their general relation is:
\begin{equation}
 \lambda_T \, \le \,  \lambda_H \, \le \,  \lambda_{GN}\, .
 \end{equation}

\noindent {\bf Remarks:}  The two first transition points were
already mentioned above.  The third has appeared in the work of Grimmett and Newman~\cite{GN}
within the context of percolation models on products of regular
trees and Euclidean lattices, where its analog  is the threshold
for the uniqueness of the infinite cluster.  The above work
motivated ref.~\cite{P}, where it was shown that
$\lambda_H<\lambda_{GN}$ for the contact process on regular trees
of degree four or more; the proof was extended to all regular
trees in~\cite{L1} and then more succinctly in~\cite{St}.

The contact process can be viewed in terms of a graphical
representation, whereby one traces the state of the infection over
the `space $\times$ time' graph,  $\G\times\R$.  Healing events
are represented by Poisson processes of intensity $1$ on the lines
of  $V\times \R$, and infection-transmission events are
represented by Poisson processes of intensity $\lambda J_{x,y}$ on
$(V\times V) \times \R$.  For the latter, the set $V\times V$
appears as the collection of directed edges, and an event at
$(e_{xy},t)$ represents a possible transmission from $x$ to $y$ at
time $t$. The set of possible sources of infection for a site $x$
at time $t$ is the set of all  points in $V\times \R$ from which
there is a path which does not backtrack in time, reaching $(x,t)$
without passing through any healing event.   We refer to this set
as $C(x,t)$. One may view it as the connected cluster of $(x,t)$
in an oriented percolation model. For brevity we denote $C\ =\
C(o,0)$.  For a more detailed description of the graphical
representation picture and its relation to the self-duality of the
contact process we refer the reader to~\cite{L2}.

The graphical representation  highlights the strong relationship
this process has with oriented percolation, and   yields the
following interpretation of the two transition points:
\begin{eqnarray}
\lambda_T &=& \inf\{\lambda: \E{|C|}=\infty\} \nonumber \\
\lambda_H &=& \inf\{\lambda: \, \P{  |C| = \infty   } > 0   \, \}
\, ,
\end{eqnarray}
where $|C|$ denotes the set's size. For a discrete set like $A_t$
the size refers to the set's cardinality, whereas for a generic
$S\subset V\times\R$, such as $C$, we denote by $|S|$ the total
length of the set's vertical segments.

By known arguments, the small-$\lambda$ phase has the following
characteristics:

\begin{proposition}\label{prop}
For any $\lambda<\lambda_T$ there exist some $c < \infty $ and
$\tau >0$ such that
\begin{equation}\label{1.5}
\mathbb{E}^{(\{o\},0) }{(|A_t|)} \,< \, c\,  e^{-t/\tau } \, .
\end{equation}
Furthermore, if $\,\sum_x J_{o,x} \, e^{+\varepsilon |x|} \, < \,
\infty  $ for some $\varepsilon >0$, then also
\begin{equation}\label{1.6}
\mathbb{P}^{(\{o\},0)}( A_t  \cap B_r^c \neq \emptyset \, \mbox{
\, for some $t\ge 0$} ) \, \le \, k \, e^{-\mu r}  \,
\end{equation}
for some $k<\infty$ and $\mu >0$, where  $B^c_r\subset V$
represents the complement of a ball of radius $r$ around the
vertex $o$.
\end{proposition}
The proof of the above proposition follows from a subadditivity
property of the contact process and can be found in several places
in the literature. For completeness we include a proof of the
above proposition in Appendix (\ref{appA}).


\subsection{Summary of the main results}

Among the key statements proven below is:

\begin{theorem}\label{thm:T=H}
  For any contact process on a
transitive graph, with a translation-invariant
infection-transmission rate $J_{x,y}$ and a constant healing rate
1,
\begin{equation}
\lambda_T \ = \ \lambda_H \, .
\end{equation}
\end{theorem}

As discussed above, for $\mathcal{G}=\mathbb{Z}^d$ the above result was
established in~\cite{BG2}.  The main method used there, based
on~\cite{M}, readily extends to transitive graphs of
subexponential growth (such graphs are amenable, though the
converse is not true).  It allows to conclude that at $\lambda<
\lambda_H$ the probability that the infection of one site will
affect another decays exponentially in the distance (and also in time). However, if $\mathcal{G}$ is
non-amenable, e.g., a regular tree, exponential decay does not yet
imply finiteness of $\chi(\lambda)$.  Nevertheless, it is not difficult to extend
the arguments given in Section 3.2 of~\cite{BG2}, which follows the
approach of~\cite{AB}, to prove the above theorem for the full
class of transitive graphs.

In~\cite{AB}, certain non-linear differential inequalities were
derived within a somewhat natural extension of the model, for  which one adds the possibility of spontaneous
infection, at the rate $h$. Added insight is derived from the
consideration of the model within the two parameter space of
$(\lambda, h)$. The original model is then recovered through the
limit  $h\to 0$. In terms of the graphical representation of the
contact process, the spontaneous infection events are represented
by a Poisson process on $V\times \R$ with density $h \, dt$.

It may be noted that while the extra parameter $h$ has a very
natural meaning for a contact process, in the original context of
percolation it has appeared as a somewhat ad-hoc auxiliary ``ghost
field", whose introduction was motivated by an analogy with the
external magnetic field of ferromagnetic Ising spin
systems~\cite{ABF}.

Keeping the terminology used in the percolation discussion, the events of
spontaneous infection (points in space $\times$ time) will be
referred to as green sites,  and their collection denoted  by $G$.
The function $\theta_+ (\lambda)$ which referred to the limiting
density of infection starting from the `all infected' state, finds
its extension to $h>0$ in the function:
\begin{equation}
\theta(\lambda,h)=\mathbb{P}(C(o,0) \cap G \neq \emptyset) \, .
\end{equation}

Following are some of the relevant properties of this extension of
the model.

\begin{lemma}\label{lem:prep}  For any contact process on a
transitive graph: \\
\noindent {\em i.\/}  At each $h>0$  there is a unique stationary
state, to which the state of the system converges for all
asymptotic initial condition ($(S_{-T},-T)$ for $ T \to \infty$)
with the infection
density given by the above function $\theta(\lambda,h)$.    \\
\noindent {\em ii.\/}  For $h>0$, the function  $\theta(\lambda,
h) $ is monotone in its arguments, continuous in $\lambda$, and
continuously differentiable  in $h$.\\
\noindent {\em iii.\/} In the limit $h \to 0+$, the function
$\theta(\lambda, h) $ yields the  quantities which were introduced
above for the  $h=0$ model as follows:  \
\begin{eqnarray}
\lim_{h \searrow 0} \, \theta(\lambda, h) & = & \theta_+(\lambda) \\
\mbox{ for $\lambda < \lambda_H$:} \qquad \lim_{h \searrow 0} \,
\frac{\partial \theta(\lambda, h)}{\partial h} & =&
\chi(\lambda)\, .
\end{eqnarray}
\end{lemma}

Since the main idea is rather standard, we relegate the proof of
the above lemma to Appendix (\ref{sec:proof}). As can be seen
there, the graphical representation provides the following useful
expressions for $\theta$ and its derivative $\chi(\lambda,h) :=
\partial \theta(\lambda, h) /
\partial h  $ is:
\begin{eqnarray}
\label{eqntheta} \theta(\lambda, h) & = & \mathbb{P}(C(o,0) \cap G
\neq \emptyset) \ = \
\E{1-e^{-h |C|} } \, ,  \\
 \nonumber \\
\mbox{ for $\lambda < \lambda_H$:} \qquad \chi(\lambda,h) & =& \E{
\, |C(o,0)|;  \, C(o,0) \cap G  = \emptyset  \, }  \ = \ \E{|C|\,
e^{-h |C|} } \, .
\end{eqnarray}

The graphical representation enables the derivation of partial
differential inequalities which, via integration through the
two-parameter space prove Thm~\ref{thm:T=H} and provide also
additional information about the behavior in the vicinity of the
critical point, which can now be commonly denoted as $\lambda_c \
:=\ \lambda_T \ = \ \lambda_H $.
\begin{theorem} \label{thm2} For any transitive graph  \\
(i) for $\lambda > \lambda_c$:
\begin{equation}\label{lambdaexp}
\theta_+(\lambda)\, \ge \,  \mbox{Const.} \, (\lambda  -
\lambda_c)^{1}
\end{equation}
(ii) at $\lambda = \lambda_c$:
\begin{equation}\label{hexponent}
\theta(\lambda_c,h)\, \ge \,  \mbox{Const.} \,  h^{1/2}   \,.
\end{equation}
\end{theorem}

The inequalities imply  bounds for the associated critical
exponents:
\begin{equation} \label{expbounds}
\beta \ge 1\, ,   \qquad \qquad \delta \ge 2\,.
\end{equation}
  As explained above,
these results are known already for both the discrete-time contact process~\cite{AB}, and the continuous-time model~\cite{BG2}.
The main novelty here is in the derivation for the continuous-time
process of the partial differential inequalities which are
discussed next.

It should be noted that the critical exponent bounds \eqref{expbounds} are saturated for the contact process
on regular  trees (of degree three or more)~\cite{W, Sch}, and also on $\mathbb{Z}^d$ when $d$ is very large or  just $d>4$ and the kernel is sufficiently `spread-out'~\cite{BW, Sa}.
The discrete-time version of this statement was proven earlier through the combination of the results of~\cite{BA,NY}.

\subsection{The key differential inequalities}

The derivation of the above results proceeds through certain non-linear partial differential inequalities (PDI).  The simplest of these is:
\begin{equation}\label{1.10}
\frac{\partial \chi}{\partial \lambda} \, \le \, |J|\, \chi^2 \, .
\end{equation}
This relation, which for percolation was presented in \cite{AN},
is basically known in the generality considered here.  It has been
noted that \eqref{1.10}   implies a critical exponent bound
($\gamma \le 1$)  which concerns the divergence rate for $\chi$ as
$\lambda \nearrow \lambda_T$:
\begin{equation}
\chi(\lambda)\, \ge \,   \frac{|J|^{-1}} { |\lambda_T -
\lambda|_{+}} \, .
\end{equation}
Next are partial differential inequalities which are similar to
the PDI which were derived in~\cite{AB} for the discrete-time
contact process,
 in the context of percolation model,  extending an earlier differential
inequality of~\cite{CC}, which has yielded a percolation analog of
\eqref{lambdaexp} with $\lambda_c$ interpreted as $\lambda_H$.

\begin{theorem}\label{thm:PDI}
For any contact process on a transitive graph, at Lebesgue almost
every $(\lambda,h)\in \R_+\times \R_+$ (due to the monotonicity of
$\theta$, the  derivatives exist in this sense):
\begin{equation} \label{eq:PDI1a}
\frac{\partial \theta}{\partial \lambda }\ \leq \ \theta |J|
\frac{\partial \theta}{\partial h}
\end{equation}
and
\begin{equation} \label{eq:PDI2a}
\theta \ \leq \ h \frac{\partial \theta}{\partial h} \ +  \left(2
\lambda^2 |J| \theta + h\lambda \right)  \frac{\partial
\theta}{\partial \lambda} \ + \ \theta^2 \, .
\end{equation}
\end{theorem}

In ref.~\cite{AB}, where the
discrete version of the above theorem was established, it was
envisioned that an extension to the continuum ought to be possible
through a limiting argument, but the result may involve some more
complicated coefficients.
Nevertheless, as is shown here the inequalities are valid in a rather simple form, which is not that different from the discrete-time version.  The proof proceeds through a
finite-volume version of the statement, given in
Theorem~\ref{thm3} in Section~\ref{sec:PDI}. Theorem~\ref{thm:PDI}
is proved in Section~\ref{sec:analysis}.

For the purpose of the derivation let us present some notions
which are of general use when working with Poisson processes.

\section{A Poisson process differentiation formula}\label{PP}

The analysis is made clearer by recognizing a general expression
for the  derivatives of the probabilities of monotone events with
respect to Poisson densities.  It forms a continuum analog of
`Russo's formula' which applies in the  discrete  setting.

\begin{definition}\label{vardef}
Let $\X $ be  a measure space, and $\Psi(\rho)$ a monotone
functional on the space ${\mathcal M}$ of non-negative measures
$\rho(x)$ on  $\X $.   A function $K(x,\rho)$ on $\X \times
{\mathcal M}$ is said to be the variational derivative of
$\Psi(\cdot )$ at $x$ if  for all finite positive continuous
measures $\alpha$ on $\X $
\begin{equation} \label{eq:defder}
\left. \frac{d}{ds}\Psi(\rho + s\,\alpha )\right|_{s=0^+}=\int
K(x)\,d\alpha(x) \, .
\end{equation}
\end{definition}
It is easy to see that when it exists,  the variational derivative
is unique.  We denote it
\begin{equation}
\frac{\delta \Psi}{\delta\rho(x)}\ =\  K(x) \, .
\end{equation}

We shall now consider functionals of the form
\begin{equation} \label{eq:psi}
\Psi(\rho) \ = \  \mathbb{P}_{\rho}( \omega \in F)
\end{equation}
where $F$ is an increasing  event, i.e.,  one whose indicator
function is a non-decreasing function of  the configuration
$\omega\subset \X$, and the subscript on $\mathbb{P}$ indicates
that $\omega$ is distributed by the Poisson process with intensity
measure $\rho$.

\begin{definition}\label{pivotal} Let $F$ be an increasing event defined for
the point process. A point $x\in\X$ is said to be {\em pivotal for
$F$ in the configuration $\omega$} if $\omega\backslash \{ x \}
\not\in F$ but $\omega\cup \{x\}\in F$. The set of pivotal points
is denoted
\begin{equation}
\Delta F(\omega)\ :=\ \{x: \mbox{$\omega\cup\{x\} \in F $ {\rm
and} $\omega \backslash \{x\} \notin F $} \} \, .
\end{equation}
\end{definition}

\begin{lemma}\label{extrusso}
For any Poisson process, the probability of any increasing event
$F$ has a variational derivative given by
\begin{equation} \label{eq:der}
\frac{\delta\, \mathbb{P}_{\rho}(F)}{\delta \, \rho(x)}\ =\
\mathbb{P}_{\rho}( F^c ; \{x\in \Delta F  \} )  \, .
\end{equation}
If  the density $\rho$ is non-atomic ($\rho( \{x\} ) =0$ for all
$x \in \X$) then also
\begin{equation} \label{eq:der2}
\frac{\delta\, \mathbb{P}_{\rho}(F)}{\delta \, \rho(x)}\ =\
\mathbb{P}_{\rho}\left(   x \in \Delta F    \right)  \, .
\end{equation}
\end{lemma}

\begin{proof}
Since the variational addition $\alpha$ is a continuous measure, a
valid way to generate a random configuration distributed by the
Poisson process at the density $\rho_s:=\rho+s\,\alpha \, $ is to
take the union of two configurations $\omega_0$ and $\widetilde
\omega_s$, drawn independently through a pair Poisson processes at
intensities $\rho$ and $s\alpha$, correspondingly. By this
construction,
\begin{eqnarray}
\mathbb{P}_{\rho_s}(F) - \mathbb{P}_{\rho_0}(F) &=& \P{\omega_0
\notin F,\,  \omega_0 \cup \widetilde \omega_s \in F }
\nonumber \\
& = &  \  \P{ \omega_0  \notin F,\,  \omega_0 \cup \widetilde
\omega_s \in F,\,
|\widetilde \omega_s |  = 1 } \ +   \\
 \ && + \
\P{\omega_0  \notin F,\,  \omega_0 \cup \widetilde \omega_s \in F,
\, |  \widetilde \omega_s |  \ge 2 }     \nonumber
\end{eqnarray}
where the first equality is due to monotonicity of $F$, and
$|\cdot |$ denotes the cardinality of a set. The first of the two
events in the last expression coincides with the event that {\em
i.} $\omega_0  \notin F$, and  {\em ii.} $\widetilde \omega_s$ is
a one-point subset of $ \Delta F(\omega_0)$. The second term is
dominated by $[s  \, \alpha(\X)]^2$. Conditioning on $\omega_0$,
and using the explicit Poisson formula for $\widetilde \omega_s$,
one gets:
\begin{eqnarray}
\mathbb{P}_{\rho_s}(F) - \mathbb{P}_{\rho_0}(F) &=& s \, \E { F^c;
\alpha(\Delta F(\omega_0) ) \, e^{- s \alpha(\Delta F(\omega_0) )
}  } \ + \ O(s^2)
\nonumber \\
&=&  s\, \int \mathbb{P}_{\rho} (F^c; x\in \Delta F )\,  d
\alpha(x) \ + \ O(s^2) \, .
\end{eqnarray}
The first claim now readily follows.

If $\rho$ is  non-atomic then the probability that the site $x$
seen on the right in \eq{eq:der}  is occupied, and thus $F$
occurs, vanishes for each a-priori specified $x\in \X$. Hence the
condition $F^c$ can be omitted from \eq{eq:der}, which is thus
reduced to \eq{eq:der2}.
\end{proof}

One may note that some auxiliary conditions are required for an
extension of the differentiation formula \eqref{eq:defder} to
apply also to the case where $\alpha$ is not a finite measure.
E.g., for any set $F$ which is measurable at infinity the pivotal
set $\Delta F$ is a.s. empty, yet the probability of $F$ need not
be independent of $\rho$.

\section{Derivation of the Partial Differential Inequalities}
\label{sec:PDI}
\subsection{A dictionary for the contact process}

We shall now translate Lemma \ref{extrusso} to the situation at
hand. Recall that the space $\times$ time picture of the contact
process is described in terms of three  independent Poisson
processes describing the random healing events, at constant rate
$1$, the spontaneous infection events, and the random
infection-transmissions.  In discussing the partial derivatives of
the corresponding probability, we allow the latter two processes
to be inhomogeneous, i.e., of densities given by functions rather
than constants: $h_x(t)$ and $\lambda_{x,y}(t)\, J_{x,y}$. The
corresponding probability measure is denoted by
$\mathbb{P}_{\lambda,h}$.

Of particular interest will be the event $E = \{C\cap
G\neq\emptyset\}$, where $C=C(o,0)$ is the infecting cluster for a
particular site $(o,0)$ . We apply in the natural way the
terminology introduced in Definition \ref{pivotal} and say that in
a given configuration $\omega$ a site $(x,t)\in V\times \R$ is
{\it green pivotal} if a change of the green set $G$ at $(x,t)$
will have an affect on whether $\omega \in E$ or not. Likewise, we
will say that an ordered bond $e^t_{yx}$, at time t, is {\it bond
pivotal} for $E$ if the presence of a transmission-event there
will affect whether $\omega \in E$ or not.

\begin{lemma}\label{lem:A} For the contact process on an arbitrary  graph,
\begin{eqnarray}
\frac{\delta \theta(\lambda,h)}{\delta \lambda_{y,x}(t)}  &= &
J_{y,x}\ \mathbb{P}_{\lambda,h}(e^t_{yx}\text{ is bond pivotal for
}E)
\label{partial_lambda} \\  \nonumber \\
\frac{\delta \theta(\lambda,h) }{\delta h_x(t)}  &=&
\mathbb{P}_{\lambda,h}((x,t)\text{ is green pivotal for }E)
\label{theta1} \label{partial_h}  \, ,
\end{eqnarray}
and for models with $h$ constant:
\begin{equation}
h\ \frac{\partial  }{\partial h}  \theta(\lambda,h) \ = \
\mathbb{P}_{\lambda,h}(C\text{ has exactly one green site})\,.
\label{theta2} \end{equation}
\end{lemma}

\begin{proof}
The first two assertions are direct consequences of Lemma
\ref{extrusso}.

For finite graphs, equation \eqref{theta2} can be understood from
\eqref{theta1},  as we comment below.  However a  direct proof
which is not limited by the finiteness condition can be obtained
from the expression \eqref{eqntheta} for $\theta(\lambda,h)$ which
readily yields:
\begin{eqnarray}
h\  \frac{\partial  }{\partial h} \theta(\lambda,h) &=&
\E{ h|C|\, e^{-h|C|}}  \nonumber \\
&=&   \mathbb{P}_{\lambda,h}(C\text{ has exactly one green
site})\, ,
\end{eqnarray}
where the last step is an explicit Poisson process relation.
\end{proof}

\noindent{\bf Remark:\/}  It is instructive to note that equation
 \eqref{theta2} can be explained by  \eqref{theta1}  through the
 following argument.
 Let $A=\{ \omega : \, |C\cap G|=1\}$ be the event that
 $C\text{ has exactly one green site}$.  Conditioned on $A$ there is
a uniquely defined site $Y(\omega)$, for
  which the event $Y(\omega) =(x,t)$
is characterized by:
 \begin{enumerate}
\item[i.]  $(x,t)\in G$, i.e., the site is an arrival point for
 the corresponding Poisson process, \\
\item[ii.]  in the configuration $\omega$, $(x,t)$ is a pivotal
site
  for $\{C\cap G\neq\emptyset\}$, i.e., for $E$.
  \end{enumerate}

The above two statements refer to independent conditions:
 {\em i.\/} referring to the status of the site itself (or arbitrarily
 small intervals including it), and
  {\em ii.\/} expressing a property of the configuration in the
  complement of this site.
 The probability of the former event  (which  is $0$ for any
 a-priori specified $t$)  has density $h$ with respect to $dt$.
A simple approximation argument can be used to show that event
{\em ii.\/}  is asymptotically independent of
 {\em i.\/} when the uncertainty interval is shrunk to a point.
This yields the identity:
\begin{equation}\label{eq:remark}
\E{\1_A\, \delta(Y-(x,t)) } \ = \ h\, \P{(x,t)\text{ is green
pivotal for }E }\, ,
\end{equation}
which is to be interpreted in a distributional sense.

Thus,
\begin{eqnarray}
\P{|C\cap G|=1} &=& \E{\1_A  } \ = \ \sum_x \int_{-\infty}^0
\E{\1_A\, \delta(Y-(x,t)) } \, dt
\nonumber \\
&=& h\,\sum_x \int_{-\infty}^0  \P{(x,t)\text{ is green pivotal
for }E } \, dt
\label{formalsum}  \\
&=& h\,\sum_x \int_{-\infty}^0
 \left. \frac{\delta \theta(\lambda,h) }{\delta h_x(t)}
\right|_{h(\cdot)\equiv h}   \,   dt  \,  .   \nonumber
\end{eqnarray}
where the last step is by  \eqref{theta1}.  Now, in case the total
time duration of the space $\times$ time graph is finite, the last
expression yields $h \, \frac{\partial \theta}{\partial h}$, and
thus we obtain \eqref{theta2}.  In this step we are applying the
Definition~\ref{vardef}  and Lemma~\ref{extrusso} with  $\alpha$
chosen to be the Lebesgue measure $dt$.  This argument is,
however,  limited by the restriction in Lemma~\ref{extrusso} that
the variational derivative $\alpha$ be finite.


\subsection{Differential inequalities at finite cutoffs}\label{DI}

As we just saw, certain technical issues need to be addressed in
order to carry the analysis directly for an infinite graph, e.g.,
the decomposition of $\frac{\partial}{\partial h} \theta $ into
the sum which appears at the last line of \eqref{formalsum} is
valid only if $\sum_x \int 1 dt < \infty$. We shall circumvent
this problem through finite time and space cutoffs.

In order to apply arguments like the one seen above, we let
$\theta_{T,L}$ the probability that infection is present at the
origin $o$  at time $0$ due to a spontaneous infection event which
has occurred within the finite time interval $(-T,0]$ at some site
within $V_L := \{ x\in V \, : \, |x|\le L\}$. Equivalently,
$\theta_{T,L}$ is the infection probability at $(o,0)$ for the
finite subgraph $\mathcal{G}_L$ with the vertex set $V_L $, in the
state which results from $A_{-T}=\emptyset$. More generally, the
infection probability in this state at $(x,t)$ is denoted by
$\theta_{T,L}(x,t)$, and we let $\theta_{T,L}^{max} :=\max_{x\in
V_L} \theta_{T,L}(x,0)$.

Due to the abundance of parameters,  the dependence of the above
quantities on $(\lambda,h)$ will occasionally be suppressed in the
notation. As a step towards Theorem~\ref{thm:PDI} we first derive
the following finite-volume version.

\begin{theorem}\label{thm3}  On the finite graph,
$\mathcal{G}_L\times [-T,0]$, the infection density introduced
above $\theta_{T,L}\equiv \theta_{T,L}(\lambda,h)$, satisfies for
$\lambda \ge 0$ and $h>0$:
\begin{equation} \label{eq:PDI1b}
\frac{\partial }{\partial \lambda }\theta_{T,L} \ \leq \ |J|\,
\theta_{T,L}^{max} \, \frac{\partial  }{\partial h} \theta_{T,L}
\end{equation}
and
\begin{equation} \label{eq:PDI2b}
\theta_{T,L}   \ \leq \ h \frac{\partial  }{\partial h}
\theta_{T,L} \ + \left(2 \lambda^2 |J| \theta_{T,L}^{max}   +
h\lambda \right) \frac{\partial  }{\partial \lambda} \theta_{T,L}
\ + \ [\theta_{T,L}^{max}]^2 \, ,
\end{equation}
\end{theorem}


This statement is proven in the rest of this section. We start
with the first inequality, using the dictionary provided by
Lemma~\ref{lem:A}.

\begin{proof}[Proof of \eqref{eq:PDI1b}]
Applying  equation~\eqref{partial_lambda}:
\begin{eqnarray}\label{eqn10}
 \frac{\partial }{\partial \lambda} \theta_{T,L}  &=& \sum_{y,x\in V_L}
\int_{-T}^0 \left. \frac{\delta \theta_{T,L}(\lambda,h) }{\delta
\lambda_{y,x}(t)}
\right|_{\lambda(\cdot)\equiv \lambda}   \,   dt \\
&=& \sum_{y,x\in V_L} J_{y,x} \int_{-T}^0
\mathbb{P}_{\lambda,h}(e^t_{yx}\text{ is vacant and is bond
pivotal for }E)\, dt. \nonumber
\end{eqnarray}
Spelling out the condition on the right-hand side we get:
\begin{eqnarray}\label{pivbond}
 \frac{\partial }{\partial\lambda} \theta_{T,L}
&=&\sum_{{y,x}\in V_L} J_{y,x} \int_{-T}^0
\mathbb{P}_{\lambda,h}((x,t) \in C;E^c; C(y,t)
\cap G\neq  \emptyset)\ dt \\
 \nonumber
&=&  \sum_{{y,x}\in V_L} J_{y,x}
\int_{-T}^0\mathbb{P}_{\lambda,h}(C(y,t)\cap G \neq \emptyset \ |\
(x,t)\in C;E^c)\  \, \mathbb{P}_{\lambda,h}((x,t)\in C;E^c)\ dt
\end{eqnarray}
The conditional expectation is the average of the probability of
the event that there is a connecting path from $G$ to $(y,t)$ in
the complement of the cluster of sites which are reached from
$(o,0)$, moving back in time, without visiting $(x,t)$.
Conditioning on the exact extent of this cluster, we see that the
conditional probability is dominated by
$\theta_{T,L}^{max}(\lambda,h)$.  Thus,
\begin{eqnarray}
\frac{\partial }{\partial\lambda} \theta_{T,L} (\lambda,h)
&\le&\theta_{T,L}^{max}(\lambda,h)\sum_{y,x\in V_L} J_{y,x}
\int_{-T}^0 \mathbb{P}_{\lambda,h}((x,t)\in
C;E^c)\ dt\\
&=&\theta_{T,L}^{max}(\lambda,h)\,|J|\,\sum_{x\in V_L}\int_{-T}^0
\mathbb{P}_{\lambda,h}((x,t)\in
C;E^c)\ dt\\
&=& \theta_{T,L}^{max}(\lambda,h)\  |J|\ \frac{\partial }{\partial
h} \theta_{T,L}(\lambda,h) \, ,
\end{eqnarray}
which is the statement we wanted to show.
\end{proof}

Inequality  \eqref{eq:PDI2b} is a bit more involved. Start by
breaking the event $\{C\cap G \neq \emptyset\}$ into two cases:
\begin{eqnarray*}
&& A\ =\ \{C\text{ has exactly one green site}\}\\
&&B\ =\ \{|C\cap G|\ge 2 \}.
\end{eqnarray*}
By Lemma \ref{lem:A}, we have $\mathbb{P}_{\lambda,h}({A})=h
\frac{\partial }{\partial h} \theta_{T,L}.$ To estimate the
probability of $B$, we split it further.  Let us define a ``gate''
as a space $\times$ time point which, in a given configuration,
needs to be visited by all paths which connect $G$ to $(o,0)$. It
is easy to see that:
\begin{enumerate}
  \item[(a)] The collection of gates is well ordered.
  \item[(b)] The {\it last gate} is either a {\it green site}  or a
vertex of a {\it transmission bond}, ``or'' taken in a
non-exclusive sense (though the probability that both occur is
zero.) We denote the former event as $B_s$ and the latter as
$B_b$.
\end{enumerate}

\begin{lemma} \label{lem:B_s}
\begin{equation} \label{eq:B_s}
\mathbb{P}_{\lambda,h}({B_s}) \ \le \
\frac{\theta_{T,L}^{max}}{1-\theta_{T,L}^{max}}\times
 h \  \frac{\partial } {\partial h} \theta_{T,L}
\end{equation}
\end{lemma}

\begin{proof}
As was discussed already, the last  factor in \eqref{lem:B_s}
coincides with $\mathbb{P}_{\lambda,h}(A)$.  The probabilities of
$B_s$ and $A$  will be compared here through the probability
densities for the uniquely defined `markers' for the two events.
For  $A$, that role is played by the unique green site in $C$, and
the corresponding decomposition of its probability is given by
equations \eqref{eq:remark} and \eqref{formalsum}.   For $B_s$ we
note that conditioned on it there is a unique site $W(\omega) \in
V_L \times [-T,0]$  for which the event $W(\omega)=(x,t)$ has the
following characteristics

\begin{enumerate}
\item[(W1)] $(x,t)\in   G$ \item[(W2)] $(x,t)\in C$ and there is
no green site within the cluster of sites from which $(o,0)$ can
be reached without visiting $(x,t)$ \item[(W3)]  there is a green
site connected to $(x,t)$ by a path in the complement of the above
cluster.
\end{enumerate}
The condition (W1) has an infinitesimal probabiliy, of density $h$
with respect to $dt$.   Let $\Phi_{x,t}$  denote the cluster
described in (W2) and $\theta(x,t)_{\Phi^c}$ be the probability of
the event (W3) conditioned on that cluster. We have the following
analog of equation \eqref{eq:remark}:
\begin{eqnarray}   \nonumber
\E{\1_{B_s}\, \delta (W-(x,t)) }  &=&    h\  \E{\1_{[\Phi_{x,t}
\cap G = \emptyset ]}  \
\theta(x,t)_{\Phi^c} }   \nonumber \\
& \le  &   h\  \E{\1_{[\Phi_{x,t} \cap G = \emptyset ]} \ (1-
\theta(x,t)_{\Phi^c}) } \times \frac{\theta_{max}}{1- \theta_{max}
}
\end{eqnarray}
where we used the fact that $\  0 \le \theta(x,t)_{\Phi^c} \le
\theta_{max} \le 1 $ and hence
\begin{equation}
\frac{\theta(x,t)_{\Phi^c} } {1-\theta(x,t)_{\Phi^c}} \ \le \
\frac{\theta_{max} } {1- \theta_{max} }
\end{equation}

Now, as is easily seen,
\begin{eqnarray}  \nonumber
\E{ \1_{[(x,t) \in C]}  \ \1_{[C \cap G = \emptyset ]} \ [1-
\theta(x,t)_{\Phi_W^c}] } & = &
   \mathbb{P}( (x,t) \in C ;\ C\cap G=\emptyset )  \\
   & = &
\mathbb{P}( (x,t) \, \mbox{is green pivotal for $E$} )
\end{eqnarray}

Putting it together, we get the following analog of
\eqref{formalsum}
\begin{eqnarray}
\mathbb{P}( B_s) & =&  \sum_{x\in V_L} \int_{-T}^0
\E{\1_{B_s}\, \delta (W-(x,t)) }  \, dt  \nonumber \\
& \le & h\, \sum_{x\in V_L} \int_{-T}^0 \mathbb{P}( (x,t) \,
\mbox{is green pivotal for $E$} ) \, dt
\, \times  \frac{\theta_{max}}{1- \theta_{max} }  \nonumber \\
& =& \frac{\theta_{max}}{1- \theta_{max} }   \times h\,
\frac{\partial }{\partial h} \theta_{T,L}
\end{eqnarray}
which  proves the lemma, through a comparison with
\eqref{formalsum}.
\end{proof}

Now we come to the trickiest estimate:

\begin{lemma} \label{lem:B_I}
\begin{equation} \label{eq:B_I}
\mathbb{P}_{\lambda,h}(B_b)(\lambda,h) \ \le \   \left[2 \lambda^2
|J| \theta_{T,L}^{max}(\lambda,h) + h\lambda \right]
\frac{\partial }{\partial \lambda } \theta_{T,L}(\lambda,h)
\end{equation}
\end{lemma}

\begin{proof}

As in the last proof, we shall compare $\P{B_b}$ with
$\frac{\partial}{\partial \lambda} \theta_{T,L}$   by expressing
each of the quantities as integrals, with simple bounds relating
the two integrands.   The probability of   $B_b$ would be
decomposed similarly to that of $B_s$ there, except that we shall
also integrate over the specifics of the last event occurring at
$x$ before the time $t$.

For a site $(x,t)$, let $\tau_{x,t}$ be the time of the last event
at $x$ preceding $t$, which can be either healing, spontaneous
infection, or an infection-transmission event into $x$. By
properties of the Poisson distribution, for $(x,t)$ specified:
$\P{t-\tau_{x,t} \ge u}=e^{-(1+h+\lambda |J|) u}$, and conditioned
on the value of $\tau_{x,t}$ the probabilities of the three
possibilities for the event, have the ratios $1:h:\lambda|J|$.

If $\Phi_{x,t}$ is as in the proof of the previous lemma, then let
$K_{x,t}^{(1)}$, $K_{x,t}^{(2)}$, and $K_{x,t}^{(3)}$ be the
events that:  $\{ (x,t) \in C \ \mbox{\ and \  } \Phi_{x,t}\cap
G=\emptyset \}$ and the last event at $x$ preceding $t$ is
correspondingly: healing, spontaneous infection, or an
infection-transmission event. Due to the independence of future
from the past events, we have:
\begin{equation} \label{ratios}
\P{K_{x,t}^1} \, : \, \P{K_{x,t}^2} \, : \,\P{K_{x,t}^3} \  = \
1\, : \, h \, : \, \lambda |J|
\end{equation}

Now,  if the event $B_b$ occurs there is a unique bond
$\widetilde{W}$ for which the event $\widetilde{W}=e^t_{yx}$ is
characterized by the conditions:
\begin{enumerate}
\item[($\widetilde{W}$1)]  the bond $e^t_{yx}$ is realized as an
infection-transmission event, \item[($\widetilde{W}$2)] there is
no green site which connects to $(o,0)$ without visiting $(x,t)$
\item[($\widetilde{W}$3)] the sites $(x,t)$ and $(y,t)$ are
reached by a pair of disjoint paths from distinct green sites
$g_x$ and $g_y$, both in the complement of $\Phi_{x,t}$.
\end{enumerate}
The existence of a site with the above characteristics is in fact
equivalent to the event $B_b$.  Thus $\P{B_b}$ can be written, in
a form similar to \eqref{formalsum}, as a sum of integrals of $\E{
\, \delta (\widetilde{W}-e^t_{yx} )\,  }$.   Splitting that
further according to the characteristics of the last event at $x$
preceding $t$, for which  $K_{x,t}^{(1)}$ is not an option, we
get:
\begin{equation}  \label{B_b}
\P{B_b} =   \sum_{ \stackrel{x,y\in V_L}{k=2,3} } \lambda J_{y,x}
\int_{-T}^0 \P{K_{x,t}^{(k)} } \ \P{\mbox{$(x,\tau_{x,t})$ and
$(y,t)$ are disjointly connected to $G$} \, | \,  K_{x,t}^{(k)} }\
dt  \, .
\end{equation}
In the statement that $(x,\tau_{x,t})$ and $(y,t)$ are disjointly
connected to $G$, it is possible that one of the paths has trivial
length, e. g. the event $K^{(2)}_{x,t}$.

For $k=2$ the conditional probability in the last expression
satisfies:
\begin{eqnarray}
\P{\mbox{$(x,\tau_{x,t})$ and $(y,t)$ are disjointly connected to
$G$} \, | \,  K_{x,t}^{(2)}  } \ = \  \nonumber \\
\qquad  =  \P{\mbox{$(y,t)$ is connected to $G$ by
 a path avoiding $(x,\tau_{x,t})$ } \, | \,  K_{x,t}^{(2)}  }   \\
 \qquad  =  \P{\mbox{$(y,t)$ is connected to $G$ by
 a path avoiding $(x,\tau_{x,t})$ } \, | \,  K_{x,t}^{(1)}  }  \, .
 \nonumber
 \end{eqnarray}
The first equality holds since under the condition $K_{x,t}^{(2)}
$ the site $(x,\tau_{x,t})$ is itself green, and the second
equality expresses the fact that the conditional probability is
not affected by the type of event which occurs at
$(x,\tau_{x,t})$.

For $k=3$ the condition that $(x,\tau_{x,t})$ is infected can be
met in two ways, since the site is at the end of an infection
transmitting bond.  By the van den Berg - Kesten inequality
\cite{vdBK}, which applies to independent systems, the probability
of the disjoint occurrence of two events is dominated by the
product of their separate probabilities.   Peeling off one of the
factors, and then switching the value of $k$, we obtain:
\begin{eqnarray}
\P{\mbox{$(x,\tau_{x,t})$ and $(y,t)$ are disjointly connected to
$G$} \, | \,  K_{x,t}^{(3)}  } \  \le  \ \qquad  \qquad   \qquad
\nonumber \\
\qquad  \qquad  \qquad  \qquad  \le
 \ 2 \, \theta^{max}_{T,L} \times  \P{\mbox{$(y,t)$ is connected to $G$ by
 a path avoiding $(x,\tau_{x,t})$ } \, | \,  K_{x,t}^{(1)}  }
 \end{eqnarray}
which

After the above bounds are inserted in \eqref{B_b}, we use
\eqref{ratios} to change also the value of $k$ in the factor
$\P{K_{x,t}^{(k)} }$ appearing there.  This results in:
\begin{eqnarray}
\P{B_b} &\le &    [h+2   \theta^{max}_{T,L}\,  \lambda \,
|J|]  \times  \nonumber  \\
&&  \sum_{  x,y\in V_L } \lambda J_{y,x} \int_{-T}^0
 \P{K_{x,t}^{(1)} ;  \{\mbox{$(y,t)$ is connected to $G$ by a path
avoiding $(x,\tau_{x,t})$ }\}  }  dt  \
 \nonumber  \\  \mbox{ } \ \nonumber    \\
& \le  &   [h+2 \theta^{max}_{T,L}\,  \lambda \, |J|] \sum_{
x,y\in V_L } \lambda J_{y,x} \int_{-T}^0
\P{\mbox{$e_{(x,t)}^t$ is bond pivotal for $E$ }  } \   dt \\
&= &    \  [h+2 \theta^{max}_{T,L}\,  \lambda \, |J| ] \ \lambda \
\frac{\partial}{\partial \lambda } \theta_{T,L} \, . \nonumber
\end{eqnarray}
where the last equation is by \eqref{eqn10}.
\end{proof}

\begin{proof}[Proof of \eqref{eq:PDI2b}]
Putting the above together, we have:
\begin{eqnarray}
\theta_{T,L} \ &\leq&  \mathbb{P}_{\lambda,h}({A})  \ + \
\mathbb{P}_{\lambda,h}({B_s})  \ + \ \mathbb{P}_{\lambda,h}({B_b})
\nonumber \\
&\leq&    \ h \frac{\partial \theta_{T,L}}{\partial h} \ + \ h
\frac{\theta_{T,L}^{max}}{1-\theta_{T,L}^{max}} \frac{\partial
\theta_{T,L}}{\partial h} \ + \ \left( 2 \lambda^2 |J|
\theta_{T,L}^{max}  + h\lambda \right) \frac{\partial
\theta_{T,L}}{\partial \lambda} \, .
\end{eqnarray}
Collecting the first two terms and multiplying through by
$(1-\theta_{T,L}^{max})$, one gets  \eqref{eq:PDI2b}.
\end{proof}

This concludes the proof of Theorem \ref{thm3}.

\section{Analysis: from the PDI to the critical behavior}
\label{sec:analysis}

We shall now extend the inequalities of Thm~\ref{thm3}  to the
infinite-volume (Theorem~\ref{thm:PDI}), and then explain how they
yield the main results stated in the introduction.

\begin{proof}[Proof of Theorem~\ref{thm:PDI}]
By the monotonicity of the contact process,
\begin{equation}\label{max}
\theta^{max}_{T,L}(\lambda,h) \ \le  \theta_{T,2L}(\lambda,h) \
\le \ \lim_{T',L' \to \infty} \theta_{T',L'}(\lambda,h) \ = \
\theta(\lambda,h) \, .
\end{equation}
This relation permits us to simplify (linearize) the problem of
passage to the limit, by replacing   $\theta^{max}_{T,L}$  in the
inequalities \eqref{eq:PDI1b} and \eqref{eq:PDI2b} by the limiting
function $\theta$.  We get
\begin{equation} \label{eq:PDI1c}
\frac{\partial }{\partial \lambda } \theta_{T,L}(\lambda,h) \ \leq
\ \theta(\lambda,h)\,  |J|\, \frac{\partial }{\partial h}
\theta_{T,L}(\lambda,h)
\end{equation}
and
\begin{equation} \label{eq:PDI2c}
\theta_{T,L} \ \leq \ h \frac{\partial }{\partial h} \theta_{T,L}
\ +  \left(2 \lambda^2 |J| \, \theta + h\lambda \right)  \,
\frac{\partial}{\partial \lambda} \theta_{T,L} \ + \ \theta^2 \, .
\end{equation}
The finite-volume quantities are differentiable  for $h>0$ (in
fact analytic in $(\lambda, \cdot)$).   In order to take the
limit, we shall interpret the inequalities in a weaker sense,
 as indicators of the corresponding relations for  integrals of the
quantities over $d \lambda \, dh$ against  suitable test
functions.    General arguments permit to conclude that in this
sense  the inequalities remain valid also in the limit.

More explicitely, through integration by parts \eqref{eq:PDI1c}
can be expressed as the relation of the following Stieltjes
integrals (each over $\R_+$)  with positive, compactly supported,
test functions  $g\in C_o(\R_+, \R_+)$
\begin{equation} \label{eq:PDI1d}
- \int \left[ \int  \theta_{T,L}(\lambda,h) \, d g  \right] dh \
\leq \ - \int \left[ \int  \,  |J|\,
 \theta_{T,L}(\lambda,h) \, d[g \theta]  \right] d\lambda\, ,
\end{equation}
where $dg$ on the left is a Stieltjes integral  at fixed $h$  and
$d[g\theta ]$  on the right is a Stieltjes integral at fixed
$\lambda$. By the bounded convergence theorem, as $T,L \to
\infty$,  the integrals  converge to those of the limit.   Since
the limiting function is also monotone in  its arguments
$(\lambda,h)$, the integration by parts can be reversed in the
limit.  The ultimate conclusion, allowed since the derivatives of
monotone functions are locally absolutely integrable,   is that
the limiting inequality  \eqref{eq:PDI1a} holds in the sense of a
relation holding at Lebesgue almost every $(\lambda,h)$.  A
similar argument permits to deduce
 \eqref{eq:PDI2a} from   \eqref{eq:PDI2c}, thereby proving
Theorem~\ref{thm:PDI}.
\end{proof}

The inequalities which are established in Theorem~\ref{thm:PDI}
are very close to what was proven in ref.~\cite{AB} for the
model's discrete-time version on $\Z^d$.
 From this point on, the analysis of the PDI is identical, and it is
covered by the general results of Lemma 4.1 and Lemma 5.1 of
ref.~\cite{AB}, which yield the following statement (formulated
here in the notation of ref.~\cite{AB}).

\begin{proposition}\label{thm:M}
Let $M(\beta,h) : \R^2 \mapsto \R\  $ be a positive function which
for $h=0$ is continuous from above and for $h>0$ is continuous,
increasing in each of its   arguments, and satisfies (in the  a.e.
sense):
\begin{eqnarray} \label{eq:M1}
  \ \frac{\partial M  }{\partial \beta } & \leq  & \phi  \,
M\, \frac{\partial M }{\partial h} \\
M & \leq   &  h \, \frac{\partial M }{\partial h} \ + \psi  \, M^a
\,
 \frac{\partial M}{\partial \beta} \ + \ M^2 \, ,
\label{eq:M2}   \end{eqnarray} with some $0<a < \infty$ and some
$\phi(\beta,h)$ and $\psi(\beta,h) $ which are finite on compact
subsets of $\R_{+}\times \R_{+}$. If there exists a value
$\beta_0\  $ for which
\begin{equation}
\lim_{h\searrow 0} M(\beta_0,h)/h \, = \, \infty
\end{equation}
then for $h \searrow 0$
\begin{equation} \label{eq:Mh}
M(\beta_0,h) \  \ge \  c_1 \, h^{1/(1+a) }
\end{equation}
and for $\beta \ge \beta_0$
\begin{equation}\label{eq:Mbeta}
M(\beta,0)  \  \ge \  c_2 \, |\beta -\beta_0|_+^{1/a} \, ,
 \end{equation}
with some $c_1, c_2 < \infty$.
\end{proposition}

\noindent{\bf Remark:}   Since at first glance it may appear
surprising that hard information about the critical behavior can
be obtained from ``soft '' inequalities  like \eqref{eq:M1} and
\eqref{eq:M2}, let us outline here the heuristics behind
Theorem~\ref{thm:M}.

First, combining   \eqref{eq:M1} and \eqref{eq:M2}, one gets:
\begin{equation} \label{eq:Mcombined}
M\  \leq   \   h \, \frac{\partial M }{\partial h} \ +
\phi(\beta,h)\, \psi(\beta,h) \, M^{(1+a)} \,
 \frac{\partial M}{\partial h} \ + \ M^2 \, .
  \end{equation}
We shall apply this relation to study the $h$ dependence at small
$h$ in the vicinity of $\beta_0$, which is analogous to our
$\lambda_T$.

It may be noted that   the inequality \eqref{eq:Mcombined} does
not add much information about the regime where $M$ is linear in $
h$ since there $h \frac{\partial}{\partial h} M \approx M$, and
thus already the first term on the right accounts for the left
side. However, at $\beta_0$  the dependence of $M$ on $h$ is
singular, and may be given by a power law: $M(\beta_0,h) \approx
h^{1/\delta}$, with some $\delta > 1$ ( the physicists convention
for the corresponding exponent).   For a shortcut, which is of
course not made in the actual proof, let us allow such an
assumption --  taken in the literal sense that $h
\frac{\partial}{\partial h} M \approx \frac{1}{\delta}M$. We now
see  that at $\beta_0$ and $h$ small,  \eqref{eq:Mcombined} holds
not because of  the first term on the right, but due to the
presence of the second:
\begin{equation} \label{eq:M3}
(1-\frac{1}{\delta}) M\   \leq   \ \phi(\beta,h)\, \psi(\beta,h)
\, M^{(1+a)} \,
 \frac{\partial M}{\partial h} \ + \ o(M) \, ,
  \end{equation}
Dividing by $M$ and integrating from $h=0$ up, one gets
\eqref{eq:Mh}, which leads to the interesting conclusion that
there is a gap in the allowed values of the exponent by which $M$
may vanish: it either vanishes linearly in $h$ or at a slower
power, $1/\delta$, with $\delta \ge 1+a$ (in our case $a=1$).

Once it is known that for $\beta \ge \beta_0$: $M(\beta,h) \ge c
\, h^{1/(1+a)}$,   a similar treatment of \eqref{eq:M2} yields
for that regime
\begin{equation}
(1-\frac{1}{1+a}) M\   \leq   \
  \psi(\beta,0) \, M \,
 \frac{\partial M}{\partial \beta} \ + \ o(M) \, .
  \end{equation}
Dividing by $M$, and integrating from $\beta_0$ upward, one gets
\eqref{eq:Mbeta}.  In particular, one learns that \\
$M(\beta,0+) >0$ for any $\beta> \beta_0\ $!

The complete proof of Theorem~\ref{thm:M}, which does not rely on
the power law assumption, can be obtained through the integration
of the inequalities \eqref{eq:M1} and \eqref{eq:M2} through
suitable regimes in the $(\beta,h)$ plane, as is done in Lemmas
4.1 and 5.1 of \cite{AB}.

\noindent
\begin{proof}[{\bf Proofs of Theorem~\ref{thm:T=H} and
Theorem~\ref{thm2}}] The two statements follow now by applying the
principle expressed in Theorem~\ref{thm:M}, to the inequalities of
Theorem~\ref{thm:PDI},  with
the correspondence: \\
$(\theta, \lambda, h) \mapsto (M,\beta,h)$. For this purpose we
note that  by a simple estimate of the contact process on a graph
with only one vertex, $\theta(\lambda,h) \ge h(1+h)$ or
equivalently $h \le \theta / (1-\theta)$, and hence inequality
\eqref{eq:PDI2a} can be brought to the form \eqref{eq:M2} with
$\psi  = 2 \lambda^2 |J| + \lambda/(1-\theta)$ and $a=1$.
\end{proof}

\section{Remarks}
\noindent {\em 1.\/}  The results presented here can be extended also
to graphs which are only quasi-transitive in the following sense.
The analysis can be adapted as long as it can be shown that for
each bounded region in the $(\lambda,h)$ plane, there are $0<
c_1\le c_2 < \infty$ such that
\begin{equation}\label{coupling}
\mathbb{P}_{\lambda,h}(C(x,t)\cap G\neq \emptyset)\  \in  [c_1
{\theta}(\lambda, h) , c_2 {\theta}(\lambda, h)  ] \,
\end{equation}
uniformly in $(x,t)$.  In particular, the conclusions of
 Theorem~\ref{thm:T=H}  hold under this `weak inhomogeneity' condition.

\noindent {\em 2.\/}  Among the cases for which Theorem~\ref{thm:T=H} applies to
are the many graphs of exponential growth which are the subject of
current research.  These include hyperbolic
tessellations, Cayley graphs of non-amenable groups, and
exponentially growing amenable graphs such as the lamplighter
group and the Diestel-Leader graph (see ~\cite{LP}  and references therein.)
An example for which the
discrete-time version of Theorem \ref{thm:T=H} was recently
applied  is the thermodynamic limit of the small-world graphs,
see~\cite{DJ}.

\noindent {\em 3.\/}   The method and results presented here apply also to unoriented
percolation models on transitive graphs similar to those
considered here,  i.e. $\mathcal{G} \times \R$  with one continuum
dimension.   Similar independence of the argument from the
presence of orientation was noted in the previous related results
on the contact process, \cite{AB,BG2}.

\noindent {\em 4.\/}
A topic which our discussion did not address is whether  in addition to the general Properties 1.-4.  it is also true, for
contact process in the generality considered here, that
the upper stationary infection density vanishes at the
critical point.   An equivalent formulation is that at $\lambda=\lambda_c$ infection from a single site will almost surely die out.   Such a statement was established for $\mathcal{G}=\mathbb{Z}^d$ in the celebrated work of
Bezuidenhout and Grimmett~\cite{BG}, and their arguments can most
likely be extended to all graphs of subexponential growth
satisfying a certain homogeneity involving block structures.
The only related results known to the authors for contact processes on graphs of exponential growth
are those of ref.~\cite{MSZ} -- where the corresponding statement is proven for regular trees,   and the corresponding statement for regular percolation on Cayley graphs of
non-amenable groups, of ref.~\cite{BLPS}.

\noindent {\em 5.\/}   Finally, we note that Theorem~\ref{thm:T=H} allows to sharpen a statement which was derived in~\cite{MSZ}.   As a step towards the proof  that
infection from a single site dies out almost surely at $\lambda_H$, it is shown there, for contact processes on tree graphs,  that
\begin{equation}
 \exp(\xi(\lambda) \, t) \le
\mathbb{E}^{(\{o\},0)}(|A_t|) \le  c \exp(\xi(\lambda) \,   t) \, .
\end{equation}
at some continuous $\xi(\lambda)$, and $c< \infty $.  It is not difficult to see that $\xi(\lambda)>0 $ for $\lambda > \lambda_H$, and by Proposition~\ref{prop}
$\xi(\lambda) <0 $ for $\lambda < \lambda_T$.   Thus,
Theorem~\ref{thm:T=H}  ($\lambda_H=\lambda_T$) allows to conclude that   $\xi(\cdot)$
actually changes sign at the transition point.

\appendix

 \addcontentsline{toc}{section}{Appendix}

\section{Exponential decay in the subcritical regime}
\label{appA} For completeness, we provide here a proof that
throughout the regime $\lambda < \lambda_T$, which is
characterized by $\chi(\lambda)<\infty$,  the probability that the
infection, if is introduced at a single site, would persist for
time $t$ and/or spread over distance $L$  decays exponentially in
$t$ and $L$.  The proof uses generally  known arguments.
\begin{proof}[Proof of Proposition \ref{prop}]
By the additivity of the contact process and the transitivity of
$\mathcal{G}$,
\begin{equation}
\mathbb{E}^{(\{o\},0)}(|A_{t+s}|)\le
\mathbb{E}^{(\{o\},0)}(|A_{t}|) \mathbb{E}^{(\{o\},0)}(|A_{s}|)\,
.
\end{equation}
Subadditivity arguments permit to conclude that
\begin{equation}\label{sm9}
\lim_{t\rightarrow\infty}\frac{1}{t}\log
\mathbb{E}^{(\{o\},0)}(|A_{t}|)=\inf_{t>0}\frac{1}{t}\log
\mathbb{E}^{(\{o\},0)}(|A_t|)=\eta
\end{equation}
exists so that $\exp(\eta t)\le \mathbb{E}^{(\{o\},0)}(|A_t|).$
Since $\lambda<\lambda_T$, it must be that $\eta<0$. For
$0<\delta<-\eta$ we can find $\bar{t}$ so that
$\mathbb{E}^{(\{o\},0)}(|A_t|)<\exp((\eta+\delta)t)$ for all
$t>\bar{t}$. Letting $1/\tau=-(\eta+\delta)$ and choosing $c$
large enough completes the proof of \eqref{1.5}.

For \eqref{1.6}, we consider a process $F^{\{o\}}_t$ which is
defined to be the contact process ignoring all healing events.  In
particular $A^{\{o\}}_s\subset F^{\{o\}}_s\subset F^{\{o\}}_t$ for
all $s<t$. We have that the left-hand side of \eqref{1.6} is
bounded by
$$\mathbb{P}^{(\{o\},0)}(x\in F_t\text{ for some }|x|>
r)+\mathbb{P}^{(\{o\},0)}(A_t\neq\emptyset).$$

Coupling $F^{\{o\}}_t$ with a branching random walk starting from
one particle at the origin gives the bound
\begin{equation}\label{lasteqn}
\mathbb{P}^{(\{o\},0)}(x\in F_t)<e^{Kt\lambda}p_t(o,x)
\end{equation} for some constant
$K$ where $p_t(o,x)$ is the transition probability of a random
walk on $\G$ (note that this is a general bound which does not
require $\lambda<\lambda_c$). A standard large deviations result
which holds for $\sum_x J_{o,x} \, e^{+\varepsilon |x|} \, < \,
\infty$, says that for all $\xi>0$ there is a $c>0$ such that
$\sum_{|x|>ut}p_t(o,x)   \le   c e^{-\xi t}$. This together with
\eqref{lasteqn} gives an exponentially decaying bound on
$\mathbb{P}^{(\{o\},0)} ( F_{r/u} \cap B_r^c \neq \emptyset )$,
whereas \eqref{1.5} implies that
$\mathbb{P}^{(\{o\},0)}(A_{r/u}\neq\emptyset)$ decays
exponentially in $r$.
\end{proof}

\section {Uniqueness of the invariant measure in the presence of
spontaneous infection} \label{sec:proof} In this appendix we prove
the basic regularity properties which were asserted for the
process at $h>0$ and its relation with the standard $h=0$ version
of the model.
\begin{proof}[Proof of Lemma~\ref{lem:prep}]
{\em i.\/ } In the presence of spontaneous infection, i.e., for
$h>0$, the introduction of initial infection at time $-T < 0$, on
a set $S$, will have negligible effect on the infection at time
$0$ as $T\to\infty$. To see this, we split the probability of
infection at a site at time $0$ into two cases, {\em (a)} when it
can be accounted for by a spontaneous infection event, and {\em
(b)} when it can be present only due to the initial conditions:
\begin{eqnarray} \nonumber
\mathbb{P}^{(B,T)}(o\in A_0) \ &=& \
\P{ C(o,0)) \cap G \cap V\times [-T,0] \neq \emptyset } +  \\
&& + \P{ C(o,0)) \cap G \cap V\times [-T,0] = \emptyset ; \,
C(o,0) \cap B\times \{-T\} \neq \emptyset }
\end{eqnarray}
The second term is negligibly small in the limit $T\to\infty$,
since
\begin{eqnarray} \nonumber
\P{ C(o,0) \cap G \cap V\times [-T,0] = \emptyset ; \,
C(o,0) \cap S\times \{-T\} \neq \emptyset }  \ =  && \\
 = \  \E{e^{-h |C(o,0) \cap G \cap V\times [-T,0]| }  ; \,
C(o,0) \cap S\times \{-T\} \neq \emptyset }  && \le  \ e^{-hT}
\label{Tago}
\end{eqnarray}
which is obtained by first conditioning on the percolation
structure, i.e., the bond variables and the healing events. Thus:
\begin{eqnarray}\label{thetafunct}
\lim_{T\to\infty} \mathbb{P}^{(B,-T)}(o\in A_0) \ &=& \
\P{ C(o,0) \cap G  \neq \emptyset } \    \\
&=& \nonumber   \E{1-e^{-h |C(o,0)|} } \  \equiv \,
\theta(\lambda, h) \, .
\end{eqnarray}
This implies assertion {\em i.\/}

{\em ii.\/} The monotonicity of $\theta(\lambda,h)$ is a standard
observation (and is valid also for the approximating functions),
and the continuous differentiability of $\theta(\lambda, h)$ in
$h$, for $h>0$ is an easy consequence \eqref{thetafunct}. We turn
our attention to the continuity of $\theta(\lambda, h)$ in
$\lambda$.

As explained above, the probability that events occurring earlier
that $T$ ago are of relevance is bounded by $e^{-hT}$. Restricting
to times $[-T,0]$: the probability that the cluster $C\cap V\times
[-T,0]$ reaches a site with $|x| \ge u T$ can in turn be bounded
by the estimates which are used in Appendix \ref{appA}. These show
that, under the assumption which is made on $\{ J_{x,y} \}$, there
exists $0<u$ such that the probability that the infection reaches
the origin from a site at distance greater than $u T$ is dominated
by $c e^{-\xi T}$. Hence, for $L=u T$:
\begin{equation}\label{uniformconv}
|\theta(\lambda,h) - \theta_{T,L}(\lambda, h)| \, \le \, e^{-hT} \
+ c e^{-\xi T} \,.
\end{equation}

The continuity of $\theta_{T,L}(\lambda, h)$ in $\lambda$ is
obvious (for a detailed argument see~\cite{L2}). Since
\eqref{uniformconv} shows that $\theta_{T,L}(\cdot, h)$ converges
uniformly to $\theta(\cdot,h)$ on $\R_+$, it must be that $\theta$
is continuous in $\lambda$.

\noindent {\em iii.\/}  The representation \eqref{thetafunct}
readily implies, via the monotone convergence theorem, that
\begin{equation}
 \theta(\lambda, 0+) \ := \  \lim_{h\searrow 0}  \theta(\lambda, h) \ = \
\P{  |C(o,0)| = \infty }  \,
\end{equation}
and for $\lambda<\lambda_H$,
\begin{equation}
\chi(\lambda)=\lim_{h \searrow 0} \, \frac{\partial
\theta(\lambda, h)}{\partial h}\ =\  \E{|C(o,0)|
1_{[|C(o,0)|\neq\infty]}}\,.
\end{equation}
 The event $\{ |C(o,0)| = \infty\}$ does not coincide with
 $\cap_{T<0} \{ o \in A_0^{(V,T)} \} $ however we claim that the
 difference is of probability zero.  More explicitly:
\begin{eqnarray}
0  & \le &   \theta(\lambda, 0+) - \theta_+(\lambda)  \  \\
& \le & \P{|C(o,0)| = \infty ,\,  \exists T\in (-\infty,0):\/
C(o,0)\cap V\times \{T\} = \emptyset }  \ = \ 0 \nonumber
\end{eqnarray}
since it is easy to see that if  $|C(o,0) \cap V\times[ T_1,0]|
=\infty $, for some $T_1<0$, then with probability $1$:
$C(o,0)\cap V\times \{T\} \neq \emptyset$ for all $T<0$.
\end{proof}
The imbedding of the contact process within the its extended
two-parameter version through the relations established in this Appendix plays a fundamental role in our analysis.   In effect, it allows to relate the regimes of $\lambda<\lambda_T$ and $\lambda>\lambda_H$, at $h=0$,  by
exploring the model along contours in the half plane $\{ (\lambda, h) \in \R_+^2 \, : \, h>0 \} $.

\end{document}